\documentclass[11pt]{article}
\usepackage{graphicx}
\usepackage{color}
\usepackage{amsmath}
\usepackage{amssymb}
\usepackage{amscd}
\usepackage{bbm}
\newcommand{\R}{\mathbb{R}}
\newcommand{\inr}[1]{\bigl< #1 \bigr>}

\newcommand{\E}{\mathbb{E}}

\newcommand{\eps}{\varepsilon}

\newcommand{\norm}[1]{\left\|#1\right\|}%

\DeclareMathOperator*{\argmin}{argmin}

\newtheorem{Theorem}{Theorem}[section]
\newtheorem{Lemma}[Theorem]{Lemma}

\newtheorem{Definition}[Theorem]{Definition}

\numberwithin{equation}{section}

\def \proof {\noindent {\bf Proof.}\ \ }

\def \endproof
{{\mbox{}\nolinebreak\hfill\rule{2mm}{2mm}\par\medbreak}}

\begin{document}

\title{{Necessary moment conditions for exact reconstruction via basis pursuit}}
\author{Guillaume Lecu\'e${}^{1,3}$  \and Shahar Mendelson${}^{2,4,5}$}

\footnotetext[1]{CNRS, CMAP, Ecole Polytechnique, 91120 Palaiseau, France.}
\footnotetext[2]{Department of Mathematics, Technion, I.I.T, Haifa
32000, Israel.}
 \footnotetext[3] {Email:
guillaume.lecue@cmap.polytechnique.fr }
\footnotetext[4] {Email:
shahar@tx.technion.ac.il}
\footnotetext[5]{Supported by the Mathematical Sciences Institute -- The Australian National University and by ISF grant 900/10.}

\maketitle

\begin{abstract}
Let $X=(x_1,...,x_n)$ be a random vector that satisfies a weak small ball property and whose coordinates $x_i$ satisfy that $\|x_i\|_{L_p} \lesssim \sqrt{p} \|x_i\|_{L_2}$ for $p \sim \log n$. In \cite{LM_compressed}, it was shown that $N$ independent copies of $X$ can be used as measurement vectors in Compressed Sensing (using the basis pursuit algorithm) to reconstruct any $d$-sparse vector with the optimal number of measurements $N\gtrsim d \log\big(e n/d\big)$. In this note we show that the result is almost optimal. We construct a random vector $X$ with iid, mean-zero, variance one coordinates that satisfies the same weak small ball property and whose coordinates satisfy that $\|x_i\|_{L_p} \lesssim \sqrt{p} \|x_i\|_{L_2}$ for $p \sim (\log n)/(\log N)$, but the basis pursuit algorithm fails to recover even $1$-sparse vectors.

The construction shows that `spiky' measurement vectors may lead to a poor performance by the basis pursuit algorithm, but on the other hand may still perform in an optimal way if one chooses a different reconstruction algorithm (like $\ell_0$-minimization). This exhibits the fact that the convex relaxation of $\ell_0$-minimization comes at a significant cost when using `spiky' measurement vectors.
\end{abstract}

In Compressed Sensing (see, e.g., \cite{MR2230846} and \cite{MR2241189}), one observes linear measurements $\inr{X_i,x_0}$, $i=1,...,N$ of an unknown vector $x_0 \in \R^n$, and the goal is to reconstruct $x_0$ {\it exactly} using those measurements, when $N$ is much smaller than $n$.

A possible recovery procedure is the basis pursuit algorithm (see, for example, \cite{MR1154788}), which is defined by
\begin{equation*}
  {\rm argmin}\big(\|t\|_1:\Gamma t=\Gamma x_0\big).
\end{equation*}
The matrix $\Gamma=N^{-1/2}\sum_{i=1}^N \inr{X_i,\cdot}e_i$ is called the measurements matrix.

The central question in this context is to construct measurements matrices $\Gamma$ for which {\it any} $d$-sparse vector $x_0$ can be reconstructed using the data $\Gamma x_0$, and using the computationally friendly basis pursuit algorithm.

\begin{Definition}\label{def:ER}
A matrix $\Gamma\in\R^{N\times n}$ satisfies the exact reconstruction property of order $d$ if for any $d-$sparse vector $x_0\in\R^n$,
  \begin{equation}
    \argmin\big(\norm{t}_1 : \Gamma t=\Gamma x_0 \big)=\{x_0\}. \tag{ER(d)}
  \end{equation}
\end{Definition}

It follows from Proposition~2.2.18 in \cite{MR3113826} that if $\Gamma$ satisfies ER($d$) then $N\gtrsim d \log\big(en/d\big)$.  Moreover, there are constructions of random matrices that satisfy ER($d$) with high probability and with the optimal number of measurements (rows) $N\sim d \log\big(en/d\big)$ (see, for example, \cite{MR2230846,MR2417886,Felix-Shahar-Holger}).

A typical example of such a matrix is the gaussian matrix, that has independent standard normal random variables as entries. However, it was not obvious whether similar reconstruction properties are true for matrices with iid entries that have heavier tails. Such matrices are `spiky' in the sense that they are likely to have several very large entries.

It should be noted that if one is interested in sparse recovery by less computationally friendly methods than the basis pursuit algorithm, one may use measurement vectors that need not have any moment. The following condition that was  recently used in \cite{Shahar-Vladimir} and \cite{Shahar-COLT} actually suffices to ensure exact reconstruction. 

\begin{Definition}\label{def:small_ball_d_sparse}
 A random vector $X$ in $\R^n$ satisfies the small ball property in the set of $d$-sparse vectors if there exist $u,\beta>0$ for which, for any $d$-sparse vector $t \in \R^n$, $P\big(|\inr{X,t}|>u \norm{t}_2\big)\geq\beta$.
\end{Definition}

The small ball property is a rather minimal assumption on the measurement vector and is satisfied in fairly general situations. For example, if one of the following simple conditions holds then $X$ satisfies the small ball property with constants that depend only on $\kappa_0$ (and on $\eps$ for the first condition):
\begin{enumerate}
\item  $X$ is isotropic (i.e. for every $t\in\R^n$, $\E\inr{X,t}^2=\norm{t}^2_2$) and for some $\eps>0$ and every $d$-sparse vector $t\in\R^n$, $\norm{\inr{X,t}}_{L_{2+\eps}}\leq \kappa_0 \norm{\inr{X,t}}_{L_2}$;
\item  $X$ is isotropic and for every $d$-sparse vector $t\in\R^n$, $\norm{\inr{X,t}}_{L_2}\leq \kappa_0 \norm{\inr{X,t}}_{L_1}$;
\item $x_1,\ldots,x_n$ are $n$ independent, real valued random variables that are absolutely continuous with respect to the Lebesgue measure and with almost surely bounded densities by $\kappa_0$ and $X=(x_1,\ldots,x_n)$. 
\end{enumerate}

For example, using Corollary~2.3 in \cite{LM_compressed}, one may show the following.

\begin{Theorem} \label{thm:ERM}There exists absolute constants $c_0,c_1$ and $c_2$ for which the following holds. Let $X$ be a random vector in $\R^n$ that satisfies the small ball property as in Definition~\ref{def:small_ball_d_sparse}, let $X_1,\ldots,X_N$ be $N$ independent copies of $X$ and set $\Gamma=\frac{1}{\sqrt{N}} \sum_{i=1}^N \inr{X_i,\cdot}e_i$. If $N \geq c_0 d \log(en/d)$ then with probability larger than $1-c_1\exp(-c_2N)$, given any $d$-sparse vector $x_0$, the only $d$-sparse vector $t$ for which $\Gamma t=\Gamma x_0$ is $x_0$ itself.
\end{Theorem}

Recall that $\ell_0$-minimization is defined by  $\min\big(\norm{t}_0:\Gamma t=\Gamma x_0\big)$, where $\norm{t}_0$ is the cardinality of the support of $t$. Theorem \ref{thm:ERM} implies that under the small ball assumption, $\ell_0$-minimization recovers  any $d$-sparse vector $x_0$ from the measurements $\Gamma x_0$ if one is given the same number of measurements as the optimal number needed for the basis pursuit algorithm.

With this observation, the behaviour of the basis pursuit algorithm when faced with data generated by a `spiky' or heavy-tailed ensembles determines whether one may this computationally friendly algorithm without paying any `price' for the convex relaxation -- $\ell_1$ minimization instead of $\ell_0$ minimization.

An indication that the basis pursuit algorithm may be used `for free' even with a slightly heavier tail behaviour than the gaussian one, has been established in \cite{LM_compressed}, where it was shown that random matrices satisfying relatively weak moment conditions also satisfy ER($d$) with the optimal number of measurements.

\begin{Theorem} \label{theo:main_LM_compressed}
   There exist absolute constants $c_0$, $c_1$ and $c_2$ and for every $\alpha \geq 1/2$ there exists a constant $c_3(\alpha)$ that depends only on $\alpha$ for which the following holds. Let $X=(x_i)_{i=1}^n$ be a random vector that satisfies
\begin{description}
\item{(a)} There are $\kappa_1,\kappa_2,w>1$ such that for every $1 \leq j \leq n$,
$\|x_j\|_{L_2} =1$, and for $p=\kappa_2\log(wn)$, $\|x_j\|_{L_p} \leq \kappa_1 p^\alpha$.

\item{(b)} The small ball property in the set of $d$-sparse vectors (as in Definition~\ref{def:small_ball_d_sparse}) is satisfied by $X$ for some $u,\beta>0$.
\end{description}
If $N \geq c_0\max\left\{d \log (en/d), (c_3(\alpha)\kappa_1)^2 (\kappa_2\log(wn))^{\max\{2\alpha-1,1\}} \right\}$ and $X_1,...,X_N$ are independent copies of $X$, then,  with probability at least $1-2\exp(-c_1\beta^2N)-1/w^{\kappa_2} n^{\kappa_2-1}$, $\Gamma=N^{-1/2}\sum_{i=1}^N \inr{X_i,\cdot}e_i$ satisfies ER($c_2u^2 \beta d$).
\end{Theorem}

For example (which may be relaxed ever further), if $x$ has mean zero, variance one and $\|x\|_{L_4} \leq \kappa \|x\|_{L_2}$, then $X=(x_1,...,x_n)$ whose coordinates are independent copies of $x$, satisfies (b) for an absolute constant $u$ and a constant $\beta$ that depends only on $\kappa$.
Hence, if $x$ is a mean-zero, variance one random variable for which $\|x\|_{L_4} \leq \kappa$ and $\|x\|_{L_p} \leq \kappa \sqrt{p}$ for $p \sim \log n$, both conditions (a) and (b) hold for $X=(x_1,...,x_n)$, and with high probability, if $N =c_1(\kappa) d \log(eN/d)$, $\Gamma$ satisfies ER($c_2(\kappa)d$).

In contrast, as noted above, the small ball assumption is the only component needed to show that  $\ell_0$-minimization  recovers $x_0$ exactly. The moment condition (a) in Theorem \ref{theo:main_LM_compressed} is used to extend the control one has from the set of $d$-sparse vectors (which is enough for the $\ell_0$-minimization procedure) to its convex hull. It is well understood that controlling the behaviour of $\Gamma$ on the convex hull of the set of $d$-sparse vectors is a key component in the analysis of the basis pursuit algorithm, and
thus, one may ask if moment properties of the measurement vector are an essential price that one has to pay to pass from $\ell_0$-minimization to its convex relaxation, the basis pursuit algorithm.

The main result of this note is a construction that shows that this is indeed the case, and that using `spiky' measurement vectors for the basis pursuit algorithm is costly.

To formulate the result, we say that a random matrix $\Gamma$ is generated by the random variable $x$ if $\Gamma=\frac{1}{\sqrt{N}}\sum_{i=1}^N \inr{X_i,\cdot}e_i$, where $X_1,...,X_N$ are independent copies of the random vector $X=(x_1,...,x_n)$ whose coordinates are independent copies of $x$.

\vskip0.5cm

\noindent {\bf Theorem A.} \textit{ There exist absolute constants $c_0,c_1,c_2$ and $c_3$ for which the following holds. Given $n\geq c_0$ and $N \log N\leq c_1 n$, there exists a mean-zero, variance one random variable $x$ that satisfies $\norm{x}_{L_4} \leq c_2$, $\norm{x}_{L_p} \leq c_2\sqrt{p}$ for $p= c_3 (\log n)/(\log N)$, and if $\Gamma$ is the $N \times n$ matrix generated by $x$ then with probability larger than $1/2$, $\Gamma$ does not satisfy the exact reconstruction property of order $1$.}

\vskip0.3cm

Note that if $\Gamma$ is generated by $x$ that satisfies $\norm{x}_{L_2}=1$, $\norm{x}_{L_4} \leq c_2$ and $\norm{x}_{L_p} \leq c_2 \sqrt{p}$ for $p\sim\log  n $, then for $N \sim \log n$,  $\Gamma$ satisfies ER($1$) with high probability. On the other hand, the random ensemble from Theorem A is generated by $x$ for which $\norm{x}_{L_2}=1$, $\norm{x}_{L_4} \leq c_2$ and $\norm{x}_{L_p} \leq c_2 \sqrt{p}$ for $p\sim (\log  n)/\log \log n $, but still does not satisfy ER($1$) with probability at least $1/2$ when $N\sim \log n$.

Therefore, in the case $d=1$, a subgaussian estimate for $p \sim \log n$ is a sharp condition for exact recovery by the basis pursuit algorithm with an optimal number of measurements (up to a $\log \log n$ factor).

\vskip0.3cm
An alternative formulation of Theorem~A is the following:
\vskip0.3cm
\noindent {\bf Theorem A$^\prime$.} \textit{If $n \geq c_0$ and $p >2$, there exists a mean-zero and variance $1$ random variable $x$, for which $\|x\|_{L_4} \leq \kappa$ and $\|x\|_{L_p} \leq \kappa\sqrt{p}$, and with probability at least $1/2$, if $N \lesssim \sqrt{p} n^{1/p}$, $\Gamma$ does not satisfy the exact reconstruction property of order 1.}
\vskip0.5cm

Observe that under the assumption of Theorem A$^\prime$, the random vector $X=(x_1,...,x_n)$ does satisfy the conditions of Theorem \ref{thm:ERM}. Therefore, one requires only $N \sim d\log(en/d)$ random measurements using independent copies of $X$ to identify any $d$-sparse vector using $\ell_0$-minimization. For $1$-sparse vectors, and, say $p=4$, the two facts imply that the price one pays for using the basis pursuit algorithm is high: $\sim n^{1/4}$ random measurements are required instead of $\sim \log n$ for the $\ell_0$ minimization.

A final remark has to do with the case in which one is given noisy measurements.
An efficient procedure in this case is the LASSO (see, e.g. \cite{MR1379242,MR2807761}). Statistical properties of the LASSO have been obtained under several hypotheses, one of which is the \textit{compatibility condition} introduced in \cite{vdG07}:
\begin{equation}
  \label{eq:compatibility_constant}
  \phi^2(L,S)=|S| \min_{\beta\in\R^n}\Big(\norm{\Gamma \beta_S-\Gamma \beta_{S^c}}_{\ell_2^N}:\norm{\beta_S}_1=1, \norm{\beta_{S^c}}_1\leq L\Big)
\end{equation}
is the compatibility constant, and the LASSO performs well when $\phi^2(L,S)$ is `large'; if $\phi^2(L,S)=0$ there are no guarantees on its performance.

The measurements matrix we will construct satisfies, with probability at least $1/2$, that $\phi^2(L,1)=0$ for any $L\geq1$. In particular, the known results on the estimation performance of LASSO under the compatibility condition (or the restricted eigenvalue assumption (cf. \cite{MR2533469})) simply do not apply.

Although we have chosen not to study the performance of the LASSO procedure in this note, it is likely that just like basis pursuit, LASSO would fail for the type of measurement matrices that we consider here, even for a $1$-sparse target vector.

\section{Proof of Theorem~A}
\label{sec:proof_theo_A}
Let $\{e_1,...,e_n\}$ be the standard basis in $\R^n$. Given an $N \times n$ matrix $\Gamma$ and $J\subset\{1,\ldots,n\}$ set $\Gamma_J$ to be the restriction of $\Gamma$ to ${\rm span}\{e_j : j \in J\}$. Let $B_1^n$ be the unit ball in $\ell_1^n$, and put $B_1^{J^c}$ to be the set of vectors in $B_1^n$ that are supported in $J^c$ -- the complement of $J$ in $\{1,...,n\}$.

\begin{Lemma}\label{lem:non_ER}
Fix integers $d,N \leq n$.
Let $v\in\R^n$ be supported on $J \subset \{1,...,n\}$ of cardinality at most $d$, that satisfies $\|v\|_1=1$. If $\Gamma v \in\Gamma B_1^{J^c}$ then  $\Gamma$ does not satisfy the exact reconstruction property of order $d$.
\end{Lemma}
\proof
Clearly, there is $w\in B_1^{J^c}$ for which $\Gamma v=\Gamma w$. Also, $v \not = w$, otherwise, $v\in B_1^J\cap B_1^{J^c}$ implying that $v=0$, which is impossible since $\norm{v}_1=1$. 

If one performs the basis pursuit algorithm trying to recover $v$ from $\Gamma v$, $w$ is at least as good `candidate' as $v$ (since $\norm{w}_1\leq 1=\norm{v}_1$), and therefore, $v$ cannot be the unique solution to the $\ell_1$-minimization problem $\min\big(\norm{t}_1:\Gamma t=\Gamma v\big)$.
\endproof

It immediately follows from Lemma~\ref{lem:non_ER} that if one wants to prove that the $N\times n$ matrix
\begin{equation*}
  \Gamma=\big(x_{ij}\big)=\left(
    \begin{array}{c}
      x_{1\cdot}^\top\\
      \vdots\\
      x_{N\cdot}^\top
    \end{array}
\right)=[x_{\cdot1},\cdots,x_{\cdot n}]
\end{equation*}
 does not satisfy ER($1$), it suffices to find $j\in\{1,\ldots,n\}$ for which
\begin{equation*}
 \Gamma e_j = x_{\cdot j} \in {\rm absconv}\big( x_{\cdot k}:k\neq j\big) = {\rm absconv}\big(\Gamma e_k: k \not = j \big).
\end{equation*}
To that end, if $B_2^N$ denotes the Euclidean unit ball in $\R^N$ and
\begin{equation}
  \label{eq:er_1_4}
  \norm{x_{\cdot j}}_2\leq \sqrt{N} \ \ \mbox{ and } \ \  \sqrt{N}B_2^N \subset {\rm absconv}\big( x_{\cdot k}:k\neq j\big),
\end{equation}
then $\Gamma$ does not satisfy ER($1$).

The proof of Theorem A and of Theorem A$^\prime$ is based on the construction of a measurements matrix for which (\ref{eq:er_1_4}) holds with probability larger than $1/2$.

Let $\eta$ be a selector (a $\{0,1\}$-valued random variable) with mean $\delta$ to be named later, and let $\eps$ be a symmetric $\{-1,1\}$-valued random variable that is independent of $\eta$. Fix $R>0$ and set $z= \eps(1+R  \eta)$.

Observe that if $p\geq2$ and $R\geq 1$ then
\begin{equation*}
  \frac{\norm{z}_{L_p}}{\norm{z}_{L_2}}= \frac{\big(1+ \big((1+R)^p-1\big)\delta\big)^{1/p}}{\big(1+ \big((1+R)^2-1\big)\delta\big)^{1/2}}\sim \frac{(1+R^p\delta)^{1/p}}{(1+R^2\delta)^{1/2}}\sim R \delta^{1/p},
\end{equation*}
provided that $R^2\delta\lesssim 1$ and that $R^p\delta\gtrsim 1$. Set $R=\sqrt{p}(1/\delta)^{1/p}$, and thus $\norm{z}_{L_p}/\norm{z}_{L_2}\sim \sqrt{p}$.

One can view $x=z/\norm{z}_{L_2}$ as a mean-zero, variance one random variable exhibiting `subgaussian' moments only up to the level $p$. Indeed, note that if $q> p$, $\norm{z}_{L_q}/\norm{z}_{L_2}\sim \sqrt{p}\delta^{1/q-1/p}$; hence, it may be far larger than $\sqrt{q}$ if $\delta$ is sufficiently small, as will be the case.

Let $X=(x_1,\ldots,x_n)$ be a vector whose coordinates are independent, distributed as $x$ and let $\Gamma$ be the measurements matrix generated by $x$.
Note that up to the normalization factor of $\norm{z}_{L_2}$, which is of the order of a constant when $R^2 \delta\lesssim 1$, $\Gamma$ is a perturbation of a Rademacher matrix by a sparse matrix with few random spikes that are either $R$ or $-R$.

Denote by $\E_\eta$ (resp. $\E_\eps$) the expectation with respect to the $\eta$-variables (resp. $\eps$-variables). A straightforward application of Khintchine's inequality (see, e.g., p.91 in \cite{LT:91}) shows that for every vector $t\in\R^n$,
\begin{align*}
  &\E\inr{X,t}^4\lesssim \E_\eta \E_\eps\Big(\sum_{j=1}^n \eps_j(1+R\eta_j)t_j\Big)^4\lesssim \E_\eta\Big(\sum_{j=1}^n (1+R\eta_j)^2t_j^2 \Big)^2
  \\
&=\E_\eta \sum_{k,\ell}(1+R\eta_k)^2t_k^2(1+R\eta_\ell)^2t_\ell^2\lesssim \norm{t}_2^4=\Big(\E\inr{X,t}^2\Big)^2
\end{align*}
provided that $R^4\delta \lesssim 1$. Applying the Paley-Zygmund theorem, it follows that the measurement vector $X$ satisfies the small ball property when  $R^4\delta\lesssim1$, and thus $\ell_0$-minimization  performs well using data generated by $X$: it requires only $cd\log(en/d)$ random measurements to reconstruct any $d$ sparse vector.

To show that the basis pursuit algorithm performs poorly using random measurements generated by $\Gamma$, set $(f_i)_{i=1}^N$ to be the canonical basis of $\R^N$ and observe that conditioned on $\eps_{ij}$'s, for every fixed $1\leq i\leq N$,
\begin{align*}
  &P_\eta\Big(\mbox{ there exists } j\in\{2,\ldots,n\}: z_{\cdot j}=\eps_{\cdot j}+\eps_{ij}R f_i \Big)\\
&=1-(1-(1-\delta)^{N-1}\delta)^{n-1}\geq 1-\frac{1}{4N}
\end{align*}
provided that
\begin{equation*}
  \frac{\log N}{n}\lesssim \delta \lesssim \frac{\log\big(en/N\big)}{N}.
\end{equation*}
Hence, by a Fubini argument, with probability at least $3/4$ there are (random) $y_1,...,y_N \in B_\infty^N$ for which
$$
{\rm absconv} \big(Rf_i+y_i : 1 \leq i \leq N\big) \subset {\rm absconv}\big( z_{\cdot k}:k\neq 1\big).
$$
\begin{Lemma} \label{lemma-width}
Using the notation above, if $v_i=Rf_i+y_{i}$ for $1\leq i\leq N$ and $y_i \in B_\infty^N$, then $\big(R/\sqrt{N}-\sqrt{N}\big) B_2^N \subset {\rm absconv}(v_1,...,v_N) \equiv V$
\end{Lemma}

\proof
A straightforward separation argument may be used to show that if, for every $w \in S^{N-1}$, $\sup_{v \in V} |\inr{v,w}| \geq \rho$, then $\rho B_2^N \subset V$ (indeed, otherwise there would be some $x \in \rho B_2^N \backslash V$; but it is impossible to separate $x$ and the convex and symmetric $V$ using any norm one functional).

Now, to complete the proof, observe that for every $w \in S^{N-1}$,
\begin{align*}
&\sup_{v \in V} |\inr{v,w}| =  \max_{1 \leq i \leq N} |\inr{Rf_i+y_{i},w}| \\
&\geq \max_{1 \leq i\leq N} |\inr{w,Rf_i}| - \max_{1 \leq i \leq N}|\inr{y_i,w}|
\geq  R/\sqrt{N} - \sqrt{N}.
\end{align*}
\endproof

Applying Lemma~\ref{lemma-width}, if $R\geq 2N$ then with probability at least $3/4$, $\sqrt{N}B_2^N\subset {\rm absconv}\big( z_{\cdot k}:k\neq 1\big)$. On the other hand, if $\delta \lesssim 1/N$ then
\begin{equation*}
  Pr[\norm{z_{\cdot1}}_2=\sqrt{N}]=(1-\delta)^N \geq 3/4.
\end{equation*}

Hence, combining the two observations, with probability at least $1/2$,
$$
\norm{z_{\cdot 1}}_2\leq \sqrt{N} \ \ {\rm and} \ \ \sqrt{N}B_2^N\subset {\rm absconv}\big( z_{\cdot k}:k\neq 1\big),
$$
and thus
\begin{equation*}
  x_{\cdot 1}\in {\rm absconv}\big(x_{\cdot k}:k\neq1).
\end{equation*}

% \texttt{[I removed all the factors ``2'' in front of $\sqrt{N}B_2^N$ because, one may just divide everything by $\norm{z}_{L_2}$ to get a result for the $x_{\cdot k}$ from the one on the $z_{\cdot k}$.]}

Of course, this assertion holds under several conditions on the parameters involved: namely, that $R=\sqrt{p}(1/\delta)^{1/p}\geq 2N$; that $(\log N)/n\lesssim \delta \lesssim\log\big(en/N\big)/N$; that $R^4\delta\lesssim 1$ and that $\delta\lesssim 1/N$.

For instance, one may select $\delta \sim (\log N)/n$ and $p\sim (\log n)/\log N$, in which case all the conditions above are met and with probability at least $1/2$, $\Gamma$ does not satisfy ER($1$), proving Theorem A. A similar calculation leads to the proof of Theorem A$^\prime$.
\endproof

\begin{footnotesize}

\bibliographystyle{plain}

\bibliography{biblio}
\end{footnotesize}

\end{document}